\documentclass[reqno, 10pt, a4paper]{amsart}

\usepackage[utf8]{inputenc}
\usepackage[OT2,T1]{fontenc}
\usepackage[english]{babel}

\usepackage{amsmath}
\usepackage{amsfonts}
\usepackage{amssymb}
\usepackage{amsthm}

\usepackage{mathrsfs}
\usepackage{listings}
\usepackage[mathcal]{euscript}
\usepackage{bbm}%

\usepackage{enumitem}
\usepackage{multirow}

\usepackage[dvipsnames]{xcolor}
\usepackage{graphicx}
\usepackage{tikz}
\usetikzlibrary{matrix}
\usepackage[all]{xy}

\usepackage{colortbl}
\definecolor{mygray}{gray}{0.92}

\usepackage{array}
\newcolumntype{C}[1]{>{\centering\arraybackslash$}p{#1}<{$}}

\usepackage{breqn}
\newcounter{myequation}[equation]

\usepackage{float}
\restylefloat{table}
\usepackage{multirow}

\usepackage{hyperref}

\usepackage{fancyvrb}

\RecustomVerbatimCommand{\VerbatimInput}{VerbatimInput}%
{
  frame=lines,
  framesep=2mm,
  fontsize=\small,
  rulecolor=\color{Gray},
  xleftmargin=1cm,
  xrightmargin=1cm
}

\theoremstyle{plain}
\newtheorem{theorem}{Theorem}[subsection]

\newtheorem{question}[theorem]{Question}

\newtheorem{proposition}[theorem]{Proposition}

\newtheorem{lemma}[theorem]{Lemma}

\newtheorem{corollary}[theorem]{Corollary}

\theoremstyle{definition}
\newtheorem{definition}[theorem]{Definition}

\newtheorem{remark}[theorem]{Remark}

\newtheorem{example}[theorem]{Example}

\def\epsilon{\varepsilon}

\def\tilde{\widetilde}

\def\Magma{\textsc{Magma}}

\DeclareMathOperator{\Aut}{Aut}

\DeclareMathOperator{\Gal}{Gal}
\DeclareMathOperator{\GL}{GL}

\DeclareMathOperator{\Proj}{Proj}

\DeclareMathOperator{\SL}{SL}

\DeclareMathOperator{\Sym}{Sym}

\newcommand{\dq}{{/\kern -3pt/}}

\def\C{\mathbb{C}}

\def\P{\mathbb{P}\,}
\def\Q{\mathbb{Q}}

\def\Z{\mathbb{Z}}

\usepackage{bm}

\def\Iv{\underline{I}}

\usepackage{cancel}

\usepackage{soul}

\hypersetup{
  pdfauthor   = {Bouchet, Lercier, Ritzenthaler, Sijsling},
  pdftitle    = {Magma functionality for isomorphisms of curves and hypersurfaces},
  pdfsubject  = {},
  pdfkeywords = {},
  colorlinks=true, breaklinks=true, urlcolor=blue, linkcolor=blue, citecolor=blue, bookmarksopen=true}

\usepackage{fullpage}
\selectfont

\usepackage[toc,page]{appendix}

\begin{document}

\title{Functionality for isomorphism classes of curves and hypersurfaces}
\date{\today}

\begin{abstract}
  We describe algorithms based on invariant theory to solve problems on the geometry of curves, mainly those of genus 2, 3 and 4. New theoretical results building on the first author's PhD thesis are also included.
\end{abstract}

\author[Bouchet]{Thomas Bouchet}
\address{%
  Thomas Bouchet,
  Center for Systems Biology Dresden,
  Max-Planck Institute of Molecular Cell Biology and Genetics,
  Pfotenhauerstrasse 108
  D-01307 Dresden,
  Germany
}
\email{bouchet@mpi-cbg.de}

\author[Lercier]{Reynald Lercier}
\address{%
  Reynald Lercier,
  DGA \& Univ Rennes, %
  CNRS, IRMAR - UMR 6625, F-35000
  Rennes, %
  France. %
}
\email{reynald.lercier@m4x.org}

\author[Sijsling]{Jeroen Sijsling}
\address{%
  Jeroen Sijsling,
  Institut für Algebra und Zahlentheorie,
  Universität Ulm,
  Helmholtzstrasse 18
  D-89081 Ulm,
  Germany
}
\email{jeroen.sijsling@uni-ulm.de}

\author[Ritzenthaler]{Christophe Ritzenthaler}
\address{%
	Christophe Ritzenthaler,
  Universit\'e C\^ote d'Azur
 Nice, %
  France. %
}
\email{christophe.ritzenthaler@univ-rennes1.fr}




\subjclass[2010]{14H10, 14H25, 14H37, 14H50}
\keywords{plane quartic curves, hyperelliptic curves, invariants, reconstruction, descent, twists}

\maketitle


This article originated in 2021 as a technical note intended to illustrate functions that complement those already available in \href{http://magma.maths.usyd.edu.au/magma/handbook/hyperelliptic_curves}{\Magma\ 2.28-27}. Most of them were disseminated across previous articles; we gathered them into coherent packages for the convenience of users. These packages, together with a full description of the corresponding functions and options, can be found at~\cite{github-hyper} for the hyperelliptic case, at~\cite{github-quart} for quartics, and at~\cite{github-g4} for curves of genus 4.

A preliminary version of this work, authored by the last three authors, was uploaded to the arXiv but never published. In the meantime, the first author completed his PhD thesis in June 2025 and obtained results that significantly extend the previously available functionality. The present article therefore provides an updated overview of the current state of these developments. Some of the new results have been published in~\cite{Bou-invs, bou-recon}, and only their algorithmic and implementation aspects are described here. However, the part of the thesis devoted to the computation of isomorphisms between hypersurfaces, and its extension to handle the case of genus 4 curves, have not yet appeared elsewhere; as it constitutes a natural continuation of the earlier version, we include the corresponding proofs in this article. Hence, this article can be read at various levels: as a vademecum on how to use the corresponding \Magma\ functions and their current limitations; as a research article with new theoretical or algorithmic results; as a source for open research questions. As a corollary, some sections may require more advanced knowledge than others and the sudden change of `pace' may be troubling for the reader not familiar with the necessary concepts of invariant theory.

Our exposition begins in Section~\ref{sec:invrec} with an overview of invariants and their role in classifying isomorphism classes of curves, both hyperelliptic and non-hyperelliptic, in genus $2$ up to $4$. Whereas the invariant theory of hyperelliptic curves was the subject of intense interest in the nineteenth century due to its relation with binary forms, methods for non-hyperelliptic curves are of more recent date, with the foundations for genus $3$ (respectively $4$) being laid in~\cite{dixmier} (respectively~\cite{Bou-invs}). After briefly touching upon the construction of these invariants, we turn to the problem of reconstruction.

The reconstruction of generic $n$-ary forms from their invariants is developed in~\cite{bou-recon}. This approach generalizes the reconstruction of  hyperelliptic curves and non-hyperelliptic genus $3$ curves from~\cite{LRS16} and applies to cases that are beyond the reach of this algorithm. It also yields a complete explicit reconstruction algorithm for generic curves of genus $4$. Finally, we present new tools for the efficient reconstruction of hyperelliptic curves of genus up to $4$.

In Section~\ref{sec:isos}, we address the explicit construction of isomorphisms between given curves. We describe new methods for computing isomorphisms between generic hypersurfaces; in particular, this yields a useful criterion for detecting when a hypersurface has trivial automorphism group (Proposition~\ref{prop:stab_trivial}). Via the correspondence with binary forms, we also obtain a new efficient method for determining isomorphisms of hyperelliptic curves. These ideas culminate in an effective algorithm for determining isomorphisms of generic hyperelliptic or non-hyperelliptic genus $4$ curves. Applications to automorphisms and twists over finite fields are discussed as well.

Finally, Appendix~\ref{sec:DOpo} contains an unpublished theoretical result from the first author which allows to certify that the reduction of the Dixmier--Ohno invariants is still a generating set of the invariant ring of ternary quartic forms in characteristic greater than 13.

The following table summarizes the state of the algorithms around rings of invariants, reconstruction, and isomorphisms, in particular the new contributions in this article. (Note that some of these are for now limited to the case of characteristic $0$, but see the list of questions below.)

\begin{center}
\footnotesize
\begin{tabular}{l|l|l|l}
& Invariants & Reconstruction & Isomorphisms \\
\hline
Hyperelliptic $g \le 4$ & Classical (\S~\ref{sec:invshyp23}) & New refinement of Mestre (\S~\ref{sec:rechyp23}) & New binary forms method (\S~\ref{sec:hypisos}) \\
Non-hyperelliptic $g = 3$ & Dixmier--Ohno (\S~\ref{sec:invsg3}) & New hypersurface method (\S~\ref{sec:recg3}) & New hypersurface method (\S~\ref{sec:isog3}) \\
Non-hyperelliptic $g = 4$ & First author \cite{Bou-invs} & First author \cite{bou-recon} & New description (\S~\ref{sec:isog4}) \\
\end{tabular}
\end{center}

\medskip
While the functionality in this article provides explicit algorithms in many practical cases, there are still several open questions remaining, as summarized in the following list.

{\bf From generic reconstruction to all curves of a given genus (already in characteristic $0$)}
\begin{itemize}
  \item In genus 3: reconstruction of all plane quartics with  automorphism group either trivial, $C_2$ or $V_4$.
  \item In genus 4: reconstruction of all curves with trivial automorphism group.
  \item In genus 4, to go further one would already need to figure out models for the different automorphisms strata of genus $4$ curves, and find the corresponding equations in the moduli space.
  \item After solving the previous problem, implement a reconstruction algorithm for curves of genus $4$ with non-trivial automorphism group.
\end{itemize}

{\bf Questions in positive characteristic}
\begin{itemize}
\item Find separants for the invariant ring of binary octic forms in characteristic $5$.
\item Prove that the separants for the invariant ring of binary octic forms in characteristic $3$ and $7$ (and $5$?) are generators.
\item Reconstruct genus $3$ hyperelliptic curves from a list of invariants (or separants) in characteristic $5$.
\item Prove that the reductions of the Dixmier--Ohno invariants are still generators for the invariant ring if the characteristic of the residue field is $11$, or $13$.
\item Prove the correctness of the conjectural homogeneous systems of parameters from~\cite{LLLR21} in characteristic $3$.
\item Determine generators for the ring of invariants of plane quartics in small characteristic.
\item Work out a reconstruction process for generic plane quartics in characteristic $2$ or $3$.
\item Find homogeneous systems of parameters in all characteristics for the ring of invariants for genus $4$ curves lying on a smooth quadric.
\item Find homogeneous systems of parameters in all characteristics for the ring of invariants for genus $4$ curves lying on a quadric cone.
\end{itemize}

\subsection*{Notation}
In this article, $k$ denotes an algebraically closed field of characteristic $p \geq 0$, whereas $K$ denotes any field.

\section{Invariants and reconstruction}\label{sec:invrec}

Since the nineteenth century, the isomorphism classes of certain algebraic varieties, notably hyperelliptic curves, have been studied through the theory of classical invariants, \textit{e.g.} binary forms under the action of the linear group. More generally, for most hypersurfaces, isomorphisms are induced by linear maps of the ambient projective space, even in positive characteristic.

\begin{itemize}
    \item Following the proof of~\cite{MatsumuraMonsky1964}, in any characteristic, isomorphisms between projective hypersurfaces of dimension at least 2  are linear, with the exception of some quartic surfaces.
    \item For plane curves, this is also the case with the exception of plane cubics using for instance~\cite[ex.18,p.56]{griffiths-harris}.
\end{itemize}

The ring of invariants of these objects under the action of the linear group is known to be finitely generated. However, explicit generators are known only in a small number of cases, especially when the characteristic of the base field is not zero.

Section~\ref{sec:invshyp23} discusses the rings of invariants of hyperelliptic curves of genus up to $4$, with some emphasis on the cases of positive characteristic. The cases of smooth non-hyperelliptic curves of genus 3 and 4 are considered in Section~\ref{sec:invsg3}.

Conversely, given an isomorphism class characterized by its invariants, one can ask for an explicit representative curve, if possible defined over the field generated by these invariants. This question is called \emph{reconstruction} (from the invariants). A general framework for addressing this question, first introduced in~\cite{bou-recon}, is summarized in Section~\ref{sec:recphil}. Concrete algorithms are considered in Section~\ref{sec:rechyp23} for hyperelliptic curves of small genus and in Section~\ref{sec:recg3} for non-hyperelliptic genus 3 and 4 curves.\smallskip

For background on invariant theory, including standard terminology and algorithms for the construction of invariants, we refer the reader to~\cite{olver, DerKem, sturmfels}. Basic definitions and properties of invariants and their generalizations (covariants and contravariants) are not recalled here.

\subsection{Invariants: hyperelliptic curves of small genus}
\label{sec:invshyp23}

We assume here that the characteristic $p$ of $k$ is different from $2$. The case $p=2$ requires separate treatment and will be addressed when possible. Let $C$ be a hyperelliptic curve of genus $g \ge 2$ over $k$ given by an affine model $y^2=f(x)$ where $f$ is a separable polynomial of degree $2g+1$ or $2g+2$.

The isomorphism class of $C$ corresponds to the orbit of the binary form $z^{2g+2} f(x/z)$ under the classical action of $\GL_2(k)$. As explained in~\cite[\S 10.2]{dolgainv}, one can characterize these orbits by the invariant space $\Proj(R_g(k))$ with
\begin{displaymath}
R_g(k) = \left( \bigoplus_{n \ge 0} \Sym^n\!\left(\Sym^{2g+2}(k^2)\right) \right)^{\SL_2(k)}\,.
\end{displaymath}

Working out generators for the invariant algebra $R_g(\C)$ was a major topic of interest among nineteenth-century mathematicians. For $g=2$, corresponding to binary sextics, such generators already appear in the work of Clebsch~\cite{clebsch}, whereas for $g=3$, that is, binary octics, they go back to~\cite{sylvester,vongall}. More recently, the case $g=4$ was settled by Brouwer and Popoviciu~\cite{brouwer2}, although substantial progress had already been made during the nineteenth century~\cite{sylvester}.

When the characteristic $p$ is positive, the situation becomes more delicate. In order to obtain a set of generators for $R_g(k)$, a good starting point is often to reduce a set of well-normalized generators of $R_g(\C)$ modulo $p$. Here well-normalized means that the invariants are primitive $\Z$-integral polynomials in the coefficients of a generic form; such a normalization is always possible by~\cite[Lemma.5.8.1]{silverman}. Unfortunately, there is currently no general criterion to ensure that the reduced invariants still generate the full invariant ring. In fact, counterexamples are known: for instance  $g=3$ and $p=5$~\cite{basson}.\smallskip

Let us give some details on the implementations available in \Magma\ for the cases $g=2$, $3$ and $4$. In all cases, the output invariants should be interpreted as points in a weighted projective space. They can be normalized and compared using the option \texttt{normalize:=true}, which relies on the techniques of~\cite[Sec.1.4]{LR11}.

\begin{itemize}

\item[] \textbf{Genus $2$.}  The algebra $R_2(\C)$ is generated by  invariants $I_2,I_4,I_6,I_{10}$ which are algebraically independent, together with an invariant $I_{15}$ such that $I_{15}^2 \in \Z[I_2,I_4,I_6,I_{10}]$. While $I_{15}$ is useful for classifying binary sextics up to  $\SL_2$-equivalence, it becomes irrelevant in the classification of genus ~2 curves, since the corresponding $\Proj R_2$ remains identical when restricting to invariants of even degree.

  Igusa~\cite{igusag2} managed to give a ``universal set of invariants'' for genus $2$ curves that works in all characteristics, including $2$. This set of invariants $\{I_2,I_4,I_6,I_8,I_{10}\}$ (with $I_8$ being superfluous beyond characteristic $2$) is integrated in \Magma\ via the function \texttt{IgusaInvariants()}\,\footnote{There are also other sets of invariants (\texttt{IgusaClebschInvariants()}, \texttt{ClebschInvariants()}) and absolute invariants (\texttt{G2Invariants()}), which are used for historical or practical reasons.}. The invariant $I_{15}$ is not computed unless the option \texttt{extend} is set to \texttt{true}.\smallskip

\item[] \textbf{Genus $3$.}
  Using the work of Geyer~\cite{geyer}, it was shown in~\cite{LR11} that the reductions modulo $p$ of Shioda generators for $R_3(\C)$ remain generators\,\footnote{Note that~\cite[Prop.1.9]{LR11} should have included a different proof of~\cite[Lemma 1]{shioda67} to be complete, as the one in \emph{loc.\ cit.}\ is only valid in characteristic $0$. The missing argument can be found in~\cite[Prop.5.5.2]{smithinvariant}.} for $R_3(k)$ whenever $p>7$. For characteristics $p = 2$, $3$ and $7$, Basson~\cite{basson} constructed sets of \emph{separating invariants} (i.e., invariants that allow one to separate the orbits of the binary forms and therefore to characterize the isomorphism classes). This list depends on the characteristic and for $p=2$ on the type defined in~\cite{NS04} or~\cite[Appendix]{basson}. For $p=5$, the algorithms return a minimal list of invariants that generates the largest subring of invariants currently known to us.

  All these variants are accessed in \Magma\ through the function \texttt{ShiodaInvariants()}. When the option \texttt{PrimaryOnly := true} is used, the output restricts to a homogeneous system of parameters (HSOP) in any characteristic different from~2 (in characteristic 2, we only conjecture that the result is a HSOP).\smallskip

\item[] \textbf{Genus $4$.}
  Brouwer and Popoviciu~\cite{brouwer2} showed that $R_4(\C)$ admits a minimal set of generators of size $106$. 
  These invariants are available in \Magma\ via the function \texttt{BrouwerPopoviciu\-Invariants()}, from the Github repository~\cite{github-g4}. To the authors' knowledge, no systematic study of the invariant ring of genus~4 hyperelliptic curves in positive characteristic was carried out so far.
\end{itemize}




\subsection{Invariants: non-hyperelliptic curves of genus 3 and 4}
\label{sec:invsg3}

Isomorphisms of non-hyperelliptic curves of genus 3 over an algebraically closed field $k$, represented by smooth plane quartics, are induced by linear transformations of the ambient projective plane $\P^2$. Accordingly, isomorphism classes are characterized by the space $\Proj(R(k))$, where
\begin{displaymath}
  R(k) = \bigoplus_{n \ge 0} \left(\Sym^n \!\left(\Sym^4(k^3)\right) \right)^{\SL_3(k)}
\end{displaymath}
is the graded ring of invariant polynomial functions on the vector space of ternary quartic forms under the classical action of $\SL_3(k)$.

When $k$ is of characteristic $0$, Dixmier~\cite{dixmier} gave a list of $7$ invariants which form a homogeneous system of parameters for $R(k)$. It was completed by Ohno~\cite{ohno}, who provided a full set of $13$ generators for the invariant ring. These invariants are polynomials in the $15$ coefficients of a ternary quartic form, with coefficients in  $\Z[1/6]$.

They can be viewed in a weighted projective space with weights $3$, $6$, $9$, $9$, $12$, $12$, $15$, $15$, $18$, $18$, $21$, $21$ and $27$. The corresponding functionality was implemented in \Magma\ for the first time in~\cite{giko}. The current function is called \texttt{DixmierOhnoInvariants()}, and differs from the earlier implementations only by normalization constants. The original normalization from~\cite{giko} can be recovered by setting the option \texttt{IntegralNormalization:=true}. Normalized representatives, for which equality suffices to test geometric isomorphism, can be obtained using the option \texttt{normalize:=true}.

\begin{remark}
  Among the Dixmier--Ohno invariants of a ternary quartic form $f(x,y,z)$, the invariant $I_{27}$ of degree $27$ plays a distinguished role. One checks that $({1}/{2^{40}})\, I_{27}$ has integral coefficients and that, over any field, its vanishing locus coincides with the locus of singular plane quartics.

  This invariant is classically computed as the resultant of the three partial derivatives of $f$, following~\cite[p.426]{gelfand}. This approach fails in characteristic $2$ for intrinsic reasons, and may have additional issues in characteristic $3$. In these cases, we instead compute the resultant of the partial derivatives of $f$ with respect to two variables and of the polynomial $f$ itself. This method is based on ideas and  programs kindly provided by Laurent Busé and relies on the techniques developed in~\cite[Def.4.6, Prop.4.7]{buse} or~\cite[Prop.11]{demazure} and~\cite[Sec.3.11.19.25]{joua}.

  This computation introduces an extraneous factor corresponding to the discriminant of the binary form $f|_{z=0}$. To eliminate this factor, we deform $f$ to $f + \epsilon (x^4+y^4)$ in characteristic different from 2 and to $f + \epsilon (x^4+xy^3+y^4)$ in characteristic $2$. The discriminant of this family is then computed and specialized at $\epsilon=0$.
\end{remark}



In his thesis~\cite{Bou-thesis}, the first author proved that, when $k$ is of characteristic $p>13$, the reductions modulo $p$ of the Dixmier--Ohno invariants remain generators of $R(k)$. The proof relies on the theory of \textit{good modules}~\cite{Don85, DM21}, together with the results of~\cite{LLLR21}, where homogeneous systems of parameters  are determined in all characteristics except~$3$ (for which only a conjectural HSOP involving an invariant of degree $81$ is proposed). The proof of this result is given in the Appendix~\ref{sec:DOpo}. We conjecture that the Dixmier--Ohno invariants also generate $R(k)$ after reduction modulo primes in the range $7< p \leq 13$. As in the case of Shioda invariants in characteristic $5$, a call to \texttt{DixmierOhnoInvariants()} in arbitrary characteristic returns a minimal set of invariants generating the largest subring of invariants currently known.
\smallskip

Invariants that classify non-hyperelliptic curves of genus 4 were worked out by the first author~\cite{Bou-invs}. Such curves lie on a quadric, and depending on whether that quadric is smooth or a cone, one gets either a set of $65$ invariants (in the smooth case) or $60$ invariants. An implementation is available in the package~\cite{github-g4}. In particular, a call to the function \texttt{InvariantsGenus4Curves()} on a smooth canonically embedded curve of genus 4 returns its associated invariants.


\subsection{Reconstruction: general philosophy}
\label{sec:recphil}


Let $V$ be a finite-dimensional vector space over an algebraically closed
field $k$, and let $G \subset \GL(V)$ be a reductive algebraic group acting
linearly on $V$.  Denote by $R$ the ring of invariant polynomial
functions on $V$.
Since invariants in $R$ enable one to classify elements in $V$ up to the
action of $G$, they allow us to distinguish
different orbits in $V$ but do not, in general, provide a canonical
representative of each orbit. This leads to the converse question, called
\emph{reconstruction from invariants}: given a point in the space of
$G$-orbits in $V$, can we find a concrete element of $V$ that corresponds to
it, up to $G$-equivalence?

The first author generalizes in his thesis earlier reconstruction algorithms
for hyperelliptic and non-hyperelliptic curves, and includes them as special cases of reconstruction for hypersurfaces.
This is explained in detail in~\cite{bou-recon}. We summarize here the general strategy and discuss its consequences for some
families of curves. It draws inspiration from Mestre's work on binary forms of
even degree~\cite{mestre}, itself based on formulae given in~\cite[\S
103]{clebsch}. The method uses the notion of contravariant, a natural
extension of invariants, and consists in constructing suitable linear bases of
contravariants combined with a Taylor-like formula, to produce a collection of
polynomials in a larger ambient space, given by the Veronese embedding of he hypersurface. The coefficients of the polynomials defining this embedding
are expressed in terms of $\mathrm{SL}_n$-invariants, and one can
generically\footnote{Precisely when the specialization of the basis of
contravariants at the set of invariants remains linearly independent.}
recover from this collection a single polynomial in the
$\mathrm{GL}_n$-equivalence class characterized by the initial invariants.
By varying the basis of contravariants, one can cover
new cases and reduce the remaining ones to small-dimensional families for
which \textit{ad hoc} methods may be applied.

When the arity and degree of the form are coprime, the reconstruction
simplifies considerably: a single polynomial is produced, which is
$\mathrm{GL}_n$-equivalent to the original one. In this situation, the
coefficients of the reconstructed polynomial are invariants, and therefore
belong to the field of moduli\footnote{There are several possible definitions
  of the field of moduli. They coincide when the base field is perfect;
  see~\cite[Sec.4]{LR11}.} of the corresponding geometric object.
This yields an explicit birational inverse to the quotient map that sends a
closed orbit to its invariants. Note that this inverse is defined on an open
set, which excludes the locus of curves that are not defined over their field
of moduli.


\begin{example}\cite[Corollary 6.4.3]{bou-recon}
\label{cor:cubic_surf}
As an illustration, we give an explicit expression for all smooth cubic
surfaces with trivial automorphism group in terms of their Clebsch-Salmon
invariants $(I_8, I_{16}, I_{24}, I_{32}, I_{40})$. The presence of
automorphisms is detected by a unique invariant of degree $100$, which also
appears (up to a constant) as the vanishing of the determinant associated with
the four generically linearly independent contravariants of order $1$.

Using the reconstruction process described above, one obtains the following polynomial:
\begin{small}
  \begin{align*}
    f = \,&(-I_8^3I_{16} + I_8^2I_{24} + 200I_{16}I_{24} + 35I_8I_{32} + 100I_{40})X^3\\[-3pt]
          & + (24I_8^2I_{16}^2 - 6I_8I_{16}I_{24} + 90I_{24}^2 - 6I_8^2I_{32} - 180I_{16}I_{32} - 30I_8I_{40})X^2Y \\[-3pt]
          & + (-12I_8^2I_{16}I_{24} + 18I_8I_{24}^2 + 240I_{24}I_{32} + 45I_8^2I_{40} - 2400I_{16}I_{40})X^2Z \\[-3pt]
          & + (-18I_8^2I_{24}^2 + 1200I_{16}I_{24}^2 + 36I_8^2I_{16}I_{32} - 24I_8I_{24}I_{32} - 720I_{32}^2 + 9I_8^3I_{40} - 360I_8I_{16}I_{40} + 960I_{24}I_{40})X^2T \\[-3pt]
          & + (-192I_8I_{16}^3 - 216I_{16}^2I_{24} - 9I_8I_{24}^2 + 66I_8I_{16}I_{32} + 60I_{24}I_{32} + 9/2I_8^2I_{40} - 60I_{16}I_{40})XY^2\\[-3pt]
          & + (192I_8I_{16}^2I_{24} + 96I_{16}I_{24}^2 - 48I_8I_{24}I_{32} - 120I_8I_{16}I_{40} - 600I_{24}I_{40})XYZ \\[-3pt]
          & + (48I_8I_{16}I_{24}^2 + 360I_{24}^3 - 576I_8I_{16}^2I_{32} - 528I_{16}I_{24}I_{32}\\[-3pt]
          & + 144I_8I_{32}^2 - 24I_8^2I_{16}I_{40} - 1920I_{16}^2I_{40} - 372I_8I_{24}I_{40})XYT \\[-3pt]
          & + (-48I_8I_{16}I_{24}^2 - 144I_{24}^3 + 480I_8I_{24}I_{40} + 1800I_{32}I_{40})XZ^2 \\[-3pt]
          & + (-144I_8I_{24}^3 + 288I_8I_{16}I_{24}I_{32} + 384I_{24}^2I_{32} + 72I_8^2I_{24}I_{40} - 8640I_{16}I_{24}I_{40} - 480I_8I_{32}I_{40} - 2400I_{40}^2)XZT \\[-3pt]
          & + (2400I_{16}I_{24}^3 + 312I_8I_{24}^2I_{32} - 432I_8I_{16}I_{32}^2 - 1056I_{24}I_{32}^2 - 2640I_8I_{16}I_{24}I_{40}\\[-3pt]
          & + 2640I_{24}^2I_{40} - 156I_8^2I_{32}I_{40} - 5280I_{16}I_{32}I_{40} - 600I_8I_{40}^2)XT^2 \\[-3pt]
          & + (512I_{16}^4 + 22I_{16}I_{24}^2 - 244I_{16}^2I_{32} + 20I_{32}^2 - 11I_8I_{16}I_{40} - 16I_{24}I_{40})Y^3 \\[-3pt]
          & + (-768I_{16}^3I_{24} - 36I_{24}^3 + 264I_{16}I_{24}I_{32} + 480I_{16}^2I_{40} + 18I_8I_{24}I_{40} - 120I_{32}I_{40})Y^2Z \\[-3pt]
          & + (-432I_{16}^2I_{24}^2 + 2304I_{16}^3I_{32} + 228I_{24}^2I_{32} - 792I_{16}I_{32}^2\\[-3pt]
          & + 216I_8I_{16}^2I_{40} - 144I_{16}I_{24}I_{40} - 114I_8I_{32}I_{40} - 120I_{40}^2)Y^2T \\[-3pt]
          & + (384I_{16}^2I_{24}^2 - 96I_{24}^2I_{32} - 240I_{16}I_{24}I_{40} + 600I_{40}^2)YZ^2 \\[-3pt]
          & + (192I_{16}I_{24}^3 - 2304I_{16}^2I_{24}I_{32} + 576I_{24}I_{32}^2 - 96I_8I_{16}I_{24}I_{40} - 1488I_{24}^2I_{40} + 1440I_{16}I_{32}I_{40} + 360I_8I_{40}^2)YZT \\[-3pt]
          & + (360I_{24}^4 - 816I_{16}I_{24}^2I_{32} + 3456I_{16}^2I_{32}^2 - 864I_{32}^3 - 2880I_{16}^2I_{24}I_{40}\\[-3pt]
          & - 360I_8I_{24}^2I_{40} + 408I_8I_{16}I_{32}I_{40} + 1344I_{24}I_{32}I_{40} + 90I_8^2I_{40}^2 - 2400I_{16}I_{40}^2)YT^2 \\[-3pt]
          & + (-64I_{16}I_{24}^3 + 400I_{24}^2I_{40} + 200I_8I_{40}^2)Z^3 \\[-3pt]
          & + (-288I_{24}^4 + 576I_{16}I_{24}^2I_{32} + 144I_8I_{24}^2I_{40} - 960I_{24}I_{32}I_{40} + 9600I_{16}I_{40}^2)Z^2T \\[-3pt]
          & + (1248I_{24}^3I_{32} - 1728I_{16}I_{24}I_{32}^2 - 10560I_{16}I_{24}^2I_{40} - 624I_8I_{24}I_{32}I_{40}\\[-3pt]
          & + 2880I_{32}^2I_{40} + 2400I_8I_{16}I_{40}^2 - 4800I_{24}I_{40}^2)ZT^2 \\[-3pt]
          & + (1600I_{16}I_{24}^4 - 1952I_{24}^2I_{32}^2 + 1728I_{16}I_{32}^3 - 1600I_8I_{16}I_{24}^2I_{40} + 2240I_{24}^3I_{40} + 1280I_{16}I_{24}I_{32}I_{40}\\[-3pt]
          & + 976I_8I_{32}^2I_{40} + 400I_8^2I_{16}I_{40}^2 - 12800I_{16}^2I_{40}^2 - 1120I_8I_{24}I_{40}^2 + 320I_{32}I_{40}^2)T^3.
  \end{align*}
\end{small}

\end{example}

\subsection{Reconstruction: hyperelliptic curves of genus 2, 3 and 4}
\label{sec:rechyp23}



In the case of hyperelliptic curves defined over a field $k$ of characteristic $p \ne 2$, the reconstruction strategy specializes to Mestre's method which we summarize here.

Given the values $\Iv$ of a set of generators for the invariants of a hyperelliptic curve $C$ of genus $g$, the reconstruction algorithm produces a conic and a plane curve of degree $g+1$, both defined over the field $K$ generated over the prime field by the invariants. Their intersection points are precisely the Weierstrass points of~$C$ and therefore determine its geometric isomorphism class.
By parameterizing the conic, possibly over a quadratic extension of~$K$, one can recover a model
of~$C$ in the usual form $y^2=f(x)$.

In order to obtain such a model defined over $K$ itself, one needs to find a $K$-rational point on the conic\,\footnote{Note that if such a point exists, then $K$ is a field of definition of $C$, but the converse  does not hold in general (see~\cite[Sec.4]{LR11}).}. In the sequel, we refer to fields for which \Magma\ can perform this step as \emph{computable fields}.

\begin{itemize}
\item[] \textbf{Genus $2$.}
For $g=2$, using the work of~\cite{CAQU,CANAPU}, reconstruction is available over any computable field,
including characteristic~2, and applies to all hyperelliptic curves of genus 2.
The corresponding functionality is implemented in \Magma\ under the name
\texttt{HyperellipticCurve\-FromIgusaInvariants()}.
\smallskip
\item[] \textbf{Genus $3$.}
Reconstruction over computable fields of characteristic $p>7$ was implemented in~\cite{LR11}
and applies to \emph{all} hyperelliptic curves of genus 3.
Thanks to the work of Basson~\cite{basson}, reconstruction from separants is also available over any
computable field of characteristic different from~5.
The corresponding function is called
\texttt{HyperellipticCurveFromShioda\-Invariants()}.
\smallskip
%
%
\item[] \textbf{Genus $4$.}
Reconstruction has been worked out by the first author for generic curves in characteristic $0$, using three generically linearly independent covariants of order $2$. The corresponding function is called \texttt{HyperellipticCurveFromBrouwerPopoviciuInvariants()}. This function works over fields of characteristic $p>6263$, although this bound could be largely improved.\smallskip

\end{itemize}



Over $\Q$, the size of the coefficients in the output of the reconstruction can be significantly reduced using methods developed by Stoll~\cite{elsenhans-stoll}, implemented in the function \texttt{MinRedBinaryForm}. Additionally, for the cases $g = 2$ and $g=3$, we speed-up the reconstruction itself.
This optimization was first introduced for $g=3$ in~\cite[Sec.~3.1]{KLLRSS17}, using a technique based on a so-called ``variation of conics'', which accelerates the search for a rational point on the conic arising in the reconstruction procedure. Developing an analogous technique for $g=2$ turns out to be more subtle, and we briefly explain the
necessary adjustments as they have not been written anywhere else.

For the usual choices of order 2 covariants used in the reconstruction,  the variation of conics involves the invariant $I_{15}$ which does not belong to the usual set of Igusa invariants.
However, there exists a relation of degree 30 that gives $I_{15}^2$ as a polynomial function of $I_2$, $I_4$, $I_6$ and $I_{10}$. If this relation is not a square over $\Q$, we can substitute $\lambda\,I_2$, $\lambda^2\,I_4$, $\lambda^3\,I_6$ and $\lambda^5\,I_{10}$ for $I_2$, $I_4$, $I_6$, \ldots $I_{10}$ for a suitable constant $\lambda$, chosen so that $\lambda^{15}\,I_{15}^2$ becomes a square. This yields $I_{15}$ up to sign, which is sufficient for our purposes. (Indeed, the genus 2 curves $y^2=f(x)$ and $y^2=-f(x)$ are quadratic twists of each other and have the same even-degree Igusa invariants, while $I_{15}(f)=-I_{15}(-f)$.)
\medskip

\begin{example}
  We reconstruct a curve of genus $3$ from its invariants.\medskip
  
  \begin{verbatim}
> P<x> := PolynomialRing(Rationals());
> f := x^8 + x^3 + 1;
> f2 := P ! (Evaluate(f, (x+1001)/(3*x+5))*(3*x+5)^8);
> f2;
6805*x^8 + 827242*x^7 + 765112882*x^6 + 306007035874*x^5 + 72331829214220*x^4 + 
    56287364680951806*x^3 + 28168431872398529278*x^2 + 8056168289692716824758*x 
    + 1008028056073190412776751
> I := ShiodaInvariants(f2);
> HyperellipticCurveFromShiodaInvariants(I);
Hyperelliptic Curve defined by y^2 = x^8 + x^3 + 1 over Rational Field
Symmetric group acting on a set of cardinality 2
Order = 2
\end{verbatim}

\end{example}
\medskip

\subsection{Reconstruction: non-hyperelliptic curves of genus 3 and 4}\label{sec:recg3}

We begin with the case of a generic non-hyperelliptic curve of genus 3 over a field $k$. The input consists of a set $\Iv$ of Dixmier--Ohno invariants. Depending on the characteristic $p$ of $k$, different instances of the general reconstruction
strategy described in Section~\ref{sec:recphil} are applied.
This leads to a generic reconstruction algorithm in characteristic different from $2$ and $3$ (see~\cite{bou-recon} for details\,\footnote{There is a typo in the definition of $p_0$ in Table~A.1 of~\cite{bou-recon}: $p_0$ should be $\mathrm{D}(C_{4,4}, c_{10,5})$, otherwise it vanishes identically.}).
\begin{itemize}

\item If $p=0$ or $p>7$ with $p \ne 17$, $37$, one uses a set of three
  contravariants of order $1$ and respective degrees $14$, $17$, and $17$. It yields an explicit plane quartic whose coefficients are rational functions in the Dixmier--Ohno invariants. Explicit expressions can be found in~\cite{github-quart}.
Over~$\Q$, the denominators of the coefficients involve only the primes $5$, $7$, $17$,
and~$37$.

\item When $p = 17$ or $37$, the same contravariants remain generically
  linearly independent, and we give alternative valid
  reconstruction formulas.

\item If $p = 5$ or $7$, one of the contravariants always vanishes. In this situation, we replace it by a contravariant of degree $20$, which again leads to valid reconstruction formulas. Using the notation of~\cite[Table A.1]{bou-recon}, this contravariant can be defined as $(\rho, \rho, c_{12,3})_2$.
\end{itemize}

When this first method fails for a given set $\Iv$, we switch to an
earlier algorithm developed in~\cite{LRS16}. This method allows one to reconstruct a model whenever the invariant $I_{12}\neq 0$.\footnote{The algorithm is only proven to work in characteristic $0$. In practice, it can be executed in positive characteristic as long as the Dixmier--Ohno invariants exist and characterize the isomorphism class, but failures may occur in characteristics as large as $89$.}

If all the above methods fail, alternative systems of covariants or contravariants can in principle
be used  to perform the reconstruction. Such systems have not been implemented, and there exist smooth plane quartics, such as the Klein quartic, for which no such system exists. Nevertheless, for all non-trivial automorphism strata except for the largest ones, namely the cyclic group $\Z/2\Z$, as well as $(\Z/2\Z)^2$ when $I_{12} = 0$, an \textit{ad hoc} reconstruction is performed and implemented using the models with few coefficients described in~\cite{LRRS}. These algorithms are implemented, although they may fail in positive characteristic as large as $41762629$.

Over $\Q$, the variation of conics together with the algorithms of~\cite{elsenhans-good} reconstruct quartics with small coefficients, as in~\cite{KLLRSS17}.
All the methods described above are implemented in the function \texttt{PlaneQuarticFromDixmierOhnoInvariants(I)}.

The case of non-hyperelliptic curves of genus~4 constitutes the main result of~\cite{bou-recon}.
At present, reconstruction is available for generic curves in characteristic $p=0$ and $p>3$, but may fail in some instances.
The corresponding functionality is implemented under the name
\texttt{Genus4CurveFromInvariants()} in the package~\cite{github-g4}.



\section{Isomorphisms, automorphisms, and twists}
\label{sec:isos}

This section deals with the computation of isomorphisms, automorphisms and twists. It is divided into three cases. We first explain a new strategy to determine (linear) isomorphisms of generic hypersurfaces. We  then give an overview of what can be achieved for hyperelliptic curves, before finishing with non-hyperelliptic curves of genus 3 and 4.

The computation of automorphisms is a special case of the one of
isomorphisms. Automorphism groups may be described as abstract groups or by
explicit linear matrices, defined either over the ground field or over a
suitable extension.

To compute twists of a quasi-projective algebraic variety $X$, one uses the bijection between the set of twists and the cohomology set $H^1(\Gal(\bar{K}/K),\Aut(X))$~\cite[p.~131]{serre-coho}. Over a finite field, this set and the bijection can be computed when the elements of $\Aut(X)$ are linear and given by explicit matrices (see for instance~\cite{MR2678623}).  This leads to a function \texttt{Twists(C, H)}, which takes as input any projective curve $C$ (not necessarily plane or non-singular), together with a finite linear subgroup $H$ of the geometric automorphism group $C$; its output yields the corresponding set of twists. We give some additional details for certain families of curves below.

\subsection{A generic strategy for hypersurfaces}\label{sec:genisos}


By further developing the idea of linear bases of covariants and contravariants introduced in  Sec.~\ref{sec:recphil}, we obtain a method to compute linear isomorphisms between generic hypersurfaces. Let $n \geq 2$, $d\geq 1$, and $\ell \geq 1$ be integers, 
where $d$ denotes the order of the  covariants that we  use, $\ell d$ is the total degree of the polynomial considered, and $n+1$ the number of variables.
Let $V$ be a $k$-vector space with basis $v_0,\ldots, v_n$, and dual basis $x_0,\ldots, x_n$. The group $\mathrm{SL}_{n+1}$ acts on $V$ by $g.v = gv$,
identifying $v$ with its coordinate vector  in the basis $(v_i)$. The induced action on $V^*$ is given by $g.x = {}^tg^{-1}x$,
where  $x$ is identified with its coordinate vector in the basis $(x_i)$.
These actions extend naturally to  $\mathrm{Sym}^d(V)$ and $\mathrm{Sym}^d(V^*)$. In what follows, we often write $\mathrm{Sym}^d(k^{n+1})$ for $\mathrm{Sym}^d(V^*)$.

Let $f,g\in\mathrm{Sym}^{\ell d}(k^{n+1})$ be the defining polynomials of two projective hypersurfaces of degree $\ell d$ and dimension $n$.
Let $c_1, \ldots, c_r$ be covariants of order $d$ and degree $\alpha$ of $\mathrm{Sym}^{\ell d}(k^{n+1})$ (we explain later how to construct them). We  assume that their evaluations at $f$ form a basis of $\mathrm{Sym}^d(k^{n+1})$; this holds for a generic $f$ when the covariants are chosen appropriately. If the same covariants evaluated at $g$ are linearly dependent,  then $f$ and $g$ are not equivalent and we are done. We now assume that this is not the case.

For every $M\in\mathrm{GL}_{n+1}(k)$,  we have
\begin{displaymath}
  c_i(M . f) = \det(M)^{\frac{\alpha \ell d - d}{n+1}} \ M . c_i(f)
  \quad \text{ for all } i.
\end{displaymath}
Let $C_f := \big(c_1(f)\ \vert \ \ldots \ \vert \ c_r(f)\big)$ be the matrix whose columns are the coordinate vectors of the $c_i(f)$ with respect to a fixed basis $\mathcal{B}$ of $\mathrm{Sym}^d(k^{n+1})$. Define $C_g$ similarly.  For  $M\in\mathrm{GL}_{n+1}(k)$,  let $\hat{M} = \rho_d(M)$, where $\rho_d: \mathrm{GL}_{n+1}(k)\longrightarrow \mathrm{GL}(\mathrm{Sym}^d(k^{n+1}))$ is the  induced representation on $\mathrm{Sym}^d(k^{n+1})$, and the matrices in the image are written with respect to the basis $\mathcal{B}$.

The following crucial lemma  immediately follows from these definitions.
\begin{lemma}\label{prop:transformation}
    If $M . f = g$, then there exists a scalar $\lambda\in k^\times$ such that $\hat{M} . C_f = \lambda C_g$.
    In particular, $\hat{M} = 1/\lambda \cdot {\ }^t(C_fC_g^{-1})$. \end{lemma}
Using Lemma~\ref{prop:transformation}, we first compute  $\hat{M}$. Recovering a preimage $M$ amounts to solving a system of monomial equations obtained from suitable entries of $\hat{M}$. We illustrate this with the following example. 

\begin{example}
    Let $n = 2$, $d = 2$, and write $M = (m_{i,j})_{1\leq i,j\leq 3}$. In the ordered basis $\mathcal{B} = \{x^2,xy,xz,y^2,yz,z^2\}$ of $\mathrm{Sym}^2(k^3)$, the induced matrix  $\hat{M} = \rho_2(M)$ is
\[\hat{M} = \begin{pmatrix}
        m_{1,1}^2 & m_{1,1}m_{1,2} & m_{1,1}m_{1,3} & m_{1,2}^2 & m_{1,2}m_{1,3} & m_{1,3}^2 \\
        2m_{1,1}m_{2,1} & \cdot & \cdot & 2m_{1,2}m_{2,2} & \cdot & 2m_{1,3}m_{2,3} \\
        2m_{1,1}m_{3,1} & \cdot & \cdot & 2m_{1,2}m_{3,2} & \cdot & 2m_{1,3}m_{3,3} \\
        m_{2,1}^2 & m_{2,1}m_{2,2} & m_{2,1}m_{2,3} & m_{2,2}^2 & m_{2,2}m_{2,3} & m_{2,3}^2 \\
        2m_{2,1}m_{3,1} & \cdot & \cdot & 2m_{2,2}m_{3,2} & \cdot & 2m_{2,3}m_{3,3} \\
        m_{3,1}^2 & m_{3,1}m_{3,2} & m_{3,1}m_{3,3} & m_{3,2}^2 & m_{3,2}m_{3,3} & m_{3,3}^2 \\
    \end{pmatrix},\]
    where the dots ``$\cdot$'' indicate non-monomial entries.

Suppose we are given the following matrix $\hat{M}_0$ in the image of $\rho_2$, \[\hat{M}_0 = \begin{pmatrix}
        0 & 0 & 0 & 5 & -5 & 5 \\
        0 & \cdot & \cdot & 20 & \cdot & -10 \\
        0 & \cdot & \cdot & 10 & \cdot & -10 \\
        5 & -10 & -5 & 20 & 10 & 5 \\
        -20 & \cdot & \cdot & 20 & \cdot & 10 \\
        20 & 10 & 10 & 5 & 5 & 5 \\
    \end{pmatrix}.\]
    To avoid introducing an unnecessary field extension, we first normalize one of the nonzero square terms to $1$, for instance the entry corresponding to $m_{1,2}^2$. We obtain
    \[\begin{pmatrix}
        0 & 0 & 0 & 1 & -1 & 1 \\
        0 & \cdot & \cdot & 4 & \cdot & -2 \\
        0 & \cdot & \cdot & 2 & \cdot & -2 \\
        1 & -2 & -1 & 4 & 2 & 1 \\
        -4 & \cdot & \cdot & 4 & \cdot & 2 \\
        4 & 2 & 2 & 1 & 1 & 1 \\
    \end{pmatrix}.\]
    Assuming $m_{1,2} = 1$ (note that $-1$ yields the same matrix), the fourth column gives $2\,m_{2,2} = 4$, $2\,m_{3,2} = 2$,
while the first row yields $m_{1,3} = -1$. The remaining monomial relations
then allow us to recover the entire matrix, leading to \[M_0 = \begin{pmatrix}
        0 & 1 & -1 \\
        -1 & 2 & 1 \\
        2 & 1 & 1
    \end{pmatrix}.\]

\end{example}
  In general, we can always choose $\mathcal{B}$ to be the standard monomial basis. The monomial entries then appear precisely in the rows and columns corresponding to pure $d$-th powers, and  the same procedure allows one to reconstruct $M$ by solving only monomial equations.

The following  proposition shows that this method applies only to hypersurfaces with  trivial automorphism group.
\begin{proposition}\label{prop:stab_trivial}
  Let $c_1, \ldots, c_r$ be covariants of $\mathrm{Sym}^{\ell d}(k^{n+1})$ of order $d$ and equal degree, whose evaluations at $f$ form a basis of $\mathrm{Sym}^d(k^{n+1})$. Then the stabilizer of $f$ in $\mathrm{PGL}_{n+1}(k)$ is trivial. Equivalently, the hypersurface $f=0$ has no nontrivial automorphisms.
\end{proposition}

\begin{proof}
    Let $M$ be in the stabilizer of $f$, \textit{i.e.} $M \cdot f = f$. By Lemma~\ref{prop:transformation}, we have $\hat{M} = \lambda {\ }^t(C_fC_f^{-1}) = \lambda {I}_r$ for some $\lambda \in k^\times$. This implies that $M$ is  diagonal, with diagonal entries $(\delta_0,\ldots, \delta_n)$. Since the diagonal entries of $\hat{M}$ are of the form $\delta_0^{d-1}\delta_i$, they must all equal $\lambda$. This forces $\delta_0 = \dots = \delta_n$, so  $M$ is a scalar multiple of the identity matrix. Therefore the stabilizer of $f$ in $\mathrm{PGL}_{n+1}(k)$ is trivial.
\end{proof}


\begin{remark}
   If the hypersurfaces defined by $f$ and $g$ have trivial automorphism groups, the  linear equivalence matrix $M$ is unique up to scalars. Its coefficients lie in the compositum of the fields of definition of $f$ and $g$, so no field extension is required to compute the isomorphism.
\end{remark}


It remains to explain how to construct such a basis of covariants. Let $m\geq 3$ be an integer, and let $d$ be the smallest divisor of $m$ for which  covariants of order $d$ of $\mathrm{Sym}^{m}(k^{n+1})$ exist, equivalently, the smallest $d$ such that $\mathrm{gcd}\left(\,{m}/{d},\ {(n+1)}/{\mathrm{gcd}(d, n+1)}\,\right) = 1$. Write $m = \ell d$.

We first fix an order $d' \leq \ell d$ (the smaller, the most convenient) satisfying the following technical conditions:
\begin{itemize}
    \item the stabilizer in $\mathrm{PGL}_{n+1}$ of a generic $f\in\mathrm{Sym}^{d'}(k^{n+1})$ is trivial; this holds for $d'>4$ when $n=1$, $d'>3$ when $n=2$ and $d'>2$ when $n>2$,
    \item $d$ divides $d'$ and
      $\mathrm{gcd}\left(\,{d'}/{d},\ {(n+1)} / {\mathrm{gcd}(d, n+1)}\,\right) = 1$.
\end{itemize}
These conditions ensure the existence of covariants of order $d$ in $\mathrm{Sym}^{d'}(k^{n+1})$ that form a basis of $\mathrm{Sym}^{d}(k^{n+1})$ generically \cite[Prop. 5.0.3]{bou-recon}.

We also require $d'$ to satisfy:
\begin{itemize}
\item the existence is constructive, providing a collection $(b_i)$ of
  covariants of order $d$ in $\mathrm{Sym}^{d'}(k^{n+1})$ that forms a basis generically. We obtain them using transvectant operations~\cite[Sec.~5]{olver} and possibly  the Clebsch transfer principle~\cite[Section~3.4.2]{dolgacag}.
\item there exists a covariant $c$ of order $d'$ for $\mathrm{Sym}^{m}(k^{n+1})$ such that, the specialization of the $(b_i)$ at $c$ forms a basis of the $k$-vector space $\mathrm{Sym}^d(k^{n+1})$.
In practice, this condition is satisfied on the first try.
\end{itemize}
This procedure  yields a collection of covariants of
$\mathrm{Sym}^{m}(k^{n+1})$ that is generically a linear basis of
$\mathrm{Sym}^d(k^{n+1})$ (see Table~\ref{tab:covorders}).
\medskip

\begin{table}[htbp]
\begin{center}
\begin{tabular}{c|c||c|c}
\textbf{Arity} & \textbf{Degree $m$} & \textbf{Order $d$'} & \textbf{Basis order $d$} \\
\hline\hline
\multirow{2}{*}{$n = 1$ (binary)}
  & $2 \mid m$    & $8$ & $2$  \\
  \cline{2-4}
  & $2 \nmid m$   & $5$ & $1$ \\
\hline
\multirow{2}{*}{$n = 2$ (ternary)}
 & $3 \mid m$     & $6$ & $3$  \\
 \cline{2-4}
 & $3 \nmid m$    & $4$ & $1$  \\
\hline
\multirow{3}{*}{$n = 3$ (quaternary)}
 & $4 \mid m$           & $4$ & $4$  \\
 \cline{2-4}
 & $2 \mid m,\ 4 \nmid m$ & $6$ & $2$  \\
 \cline{2-4}
 & $2 \nmid m$              & $3$ & $1$  \\
\hline
\end{tabular}
\end{center}
\caption{Covariant orders for generic degree-$m$ hypersurfaces in $n$
  variables.}
\label{tab:covorders}
\end{table}

\begin{remark}
    The main computational bottleneck is the symbolic computation of transvectant. Computing a level $k$ transvectant in $n+1$ variables involves a polynomial in $(n+1)^2$ variables with approximately $\binom{\left(n+1\right)^2+k-1}{k}$ terms, which grows very quickly with $n$. Faster methods for computing covariants would allow this approach to scale to higher arity.
\end{remark}

\subsection{The hyperelliptic case}
\label{sec:hypisos}

Unlike the previous section, which dealt with isomorphisms in the absence of non-trivial automorphisms, here we are interested in more global methods. As a consequence, isomorphisms are not necessarily defined over the field of definition of the curves. For both practical and theoretical reasons, we consider possibly non-algebraically closed fields $K$, call \emph{isomorphisms} those defined over $K$, and \emph{geometric isomorphisms} those defined over $k=\bar{K}$.

Throughout this section, we consider hyperelliptic curves of genus $g \ge 2$ over a base field $K$ of characteristic different from $2$, given by  explicit equations
\begin{displaymath}
  y^2 = f(x), \quad f \in K[x]\text{ separable  of degree }2g+1\text{ or }2g+2.
\end{displaymath}
Given two such curves $C_i : y^2=f_i(x)$, isomorphisms $C_1 \to C_2$ are of the form
\begin{equation} \label{eq:iso}
  (x, y) \mapsto \left(\frac{a x+b}{c x+d}, \frac{e y}{(cx+d)^{g+1}}\right),
\end{equation}
with $\left[\begin{smallmatrix} a & b \\ c & d \end{smallmatrix}\right] \in \GL_2(K)$ and $e \in K^\times$. If $K$ is algebraically closed,  one can impose $e=1$ and return to the setting of Section~\ref{sec:invshyp23}, but we do not insist on this here. The set of isomorphisms is a principal homogeneous space over the automorphism group of either curve $C_1$ or $C_2$. Determining whether $C_1$ and $C_2$ are (resp. geometrically) isomorphic, and returning the set of (resp. geometric) isomorphisms if they are, boils down to computing the elements of $\GL_2 (K)$ (resp. $\GL_2(k)$) whose right action transforms $f_1$ into a scalar multiple of $f_2$. The quotient  of the automorphism group by the hyperelliptic involution is called the group of \emph{reduced automorphisms}.

The first step is to try the fast and generic method from Section~\ref{sec:genisos}, which may succeed if the reduced geometric automorphism group of the curve is trivial. 
%
If this is not the case, another approach for determining these (geometric) isomorphisms is to apply Gröbner bases after a formal coefficient comparison, see for instance~\cite{goeb}. In~\cite{lrs}, the last three authors gave alternative algorithms that speed up this computation when $p$ does not divide $2g+2$. The function \texttt{IsIsomorphicHyperellipticCurves()} determines both the full set of isomorphisms $C_1 \to C_2$ over $K$ itself and those over the algebraic closure of $K$ using the option \texttt{geometric}. Over a number field, this option can be computationally costly. The function \texttt{AutomorphismGroup\-OfHyperellipticCurve()} returns the automorphism group of a curve as a  permutation group; if the option \texttt{explicit} is enabled, it additionally provides an isomorphism from this abstract group to the actual group of automorphisms.
\medskip

Twists of hyperelliptic curves over finite fields of characteristic different from 2  can be computed following the general strategy outlined above.
For curves of genus $2$ or $3$, the computation of twists has been thoroughly implemented in all cases, even in  characteristic $2$. For curves with small geometric automorphism groups, the generic method, recalled in the introduction of this section, is used; in the case of large automorphism groups special models are reconstructed from invariants (as in Sections~\ref{sec:invshyp23} and~\ref{sec:invsg3}), after which explicit representatives can be more easily obtained. In all cases, the corresponding function is  \texttt{Twists()}.



\subsection{The case of plane quartics}\label{sec:isog3}

Let $C_1$ and $C_2$ be two plane quartic curves over a field $K$.

As for hyperelliptic curves, the algorithm first tries to find isomorphisms $C_1 \to C_2$ using the ideas of Section~\ref{sec:genisos}. This method is quite efficient, regardless of the base field chosen, but it applies only when the chosen covariants are linearly independent, which in particular excludes curves with non-trivial automorphism groups.
Should this fail, we then check whether $I_{12} \neq 0$. In this case,~\cite{VanRijnswou} shows that the use of a suitable covariant reduces the computation of isomorphisms to the problem of finding transformations between certain binary forms associated to $C_1$ and $C_2$. This leads to the same computation of elements of $\GL_2(K)$ that was considered in Section~\ref{sec:hypisos}.

In the remaining cases, we use a direct Gröbner basis method first implemented by Michael Stoll, in which one transforms the homogeneous defining equation $F_1$ of $C_1$ by a matrix with formal coefficients and imposes that the result be equal to a nonzero scalar multiple of the defining equation $F_2$ of $C_2$. Note that, in principle, this method could be sped up by considering this equality together with those for certain covariants of $F_1$ and $F_2$ of small order. The corresponding function is called \texttt{AutomorphismGroupOfPlaneQuartic()}.

Once again, the algorithm admits both a version over the base field and a geometric version, with the latter finding the isomorphisms over the algebraic closure of $K$. Both versions are very efficient over finite fields, and the version over the base field is also reasonably fast for $K = \Q$. By contrast, finding geometric isomorphisms between plane quartic curves over $\Q$ can still take a significant amount of time. For more general fields, the implementation remains too slow, and our functions therefore restrict consideration to the cases where $K$ is either finite or equal to $\Q$.\medskip




Note that some of these routines may overlap with routines included natively in \Magma. As in the hyperelliptic case, the new ones are generally much faster.

\begin{example}
  We compare timings for our automorphism routine and the native one included in version 2.29-3 of \Magma.\medskip

  \begin{verbatim}
> P<x,y,z> := ProjectiveSpace(Rationals(),2);
> C := Curve(P, x^4+y^4+z^4);
> time aut, phi := AutomorphismGroupOfPlaneQuartic(C : geometric:=true, explicit
:= true);
Time: 0.140
> GroupName(aut);
C4^2:C3:C2
> // to get all automorphisms, we base change to Q(zeta_8)
> K := CyclotomicField(8);
> C1 := BaseChange(C,K);
> time G := AutomorphismGroup(C1);
Time: 3.200
> Gp,rep := MatrixRepresentation(G);
> GroupName(Gp);
C4^2:C3:C2
 \end{verbatim}
\end{example}

As for twists, we have implemented a function \texttt{Twists()} that computes a list of representatives of all twists of a smooth plane quartic over a finite field using the general method presented in the introduction of this section.


\subsection{The case of non-hyperelliptic curves of genus 4}
\label{sec:isog4}

We assume here that the characteristic of $K$ is different from $2$, $3$
and $5$. Non-hyperelliptic curves of genus $4$ are defined as the
intersection of a quadric $Q$ and a cubic surface $E$. Although they are not
hypersurfaces, a covariant-based approach similar to
the plane quartic case can be applied.

In the generic situation, the quadric surface $Q$ is of rank $4$. In this case, there is an extra ambiguity: $E$ can be modified by adding a linear form times $Q$, in addition to the classical action of $\GL_4(K)$. The following lemma and proposition allow us to eliminate this extra action and reduce to the standard situation without extending the field of definition.

Let $V = K^4$ with basis $\{v_0,v_1,v_2,v_3\}$ and dual basis $\{x_0,x_1,x_2,x_3\}$, equipped with the standard action of $\GL_4(K)$ as in Section~\ref{sec:genisos}.

\begin{lemma}\label{prop:trick_rank4} 
  Let $Q$ and $E$ be homogeneous polynomials of degrees $2$ and $3$ over a field $K$, with $\operatorname{ rank }(Q) =4$. Define
  \begin{displaymath}
    \tilde{E} = (Q, Q, Q, Q)_2\, E - \frac{3}{2}\, (Q, Q, Q,E)_2\, Q,
  \end{displaymath}
  where $(\cdot, \cdot, \cdot, \cdot)_l$ denotes the classical transvectant of level $l$~\cite[Section 5]{olver}. Then
  \begin{displaymath}
    \tilde{\phantom{a}} :\ (Q,E)
    \ \longmapsto\ (Q,\tilde{E})\in\mathrm{Sym}^2(V^*)\times\mathrm{Sym}^3(V^*)
  \end{displaymath}
  is  $\GL_4(K)$-equivariant, and for any linear form $\ell$ we have
  \begin{displaymath}
    \widetilde{E + \ell Q} = \tilde{E}.
  \end{displaymath}
\end{lemma}

\begin{proof}
    A direct computation shows $(Q, Q, Q, Q)_2\, \ell=({3}/{2})\,(Q, Q, Q,\ell
    Q)_2$. Equivariance of all operators then implies the lemma.
\end{proof}

\begin{proposition} 
\label{prop:isonorm}
    Let $C$ and $C'$ be two genus $4$ curves over $K$ defined by $(Q,E)$ and $(Q', E')$, respectively. We assume that $Q$ and $Q'$ are of rank $4$. Then
    \begin{displaymath}
      C\simeq C'\iff \exists\, M \in \mathrm{GL}_4(K),\ \exists\, \alpha\in K^\times,
      \quad \left\{
        \begin{array}{l}
          M . Q \,=\, Q'\,,\\
          M . \tilde{E} \,=\, \alpha \tilde{E'}\,.\\
        \end{array}
      \right.
    \end{displaymath}
\end{proposition}

\begin{proof}
    By definition, $C\simeq C'$ if and only if there exist $M\in \GL_4(K)$, $\alpha\in K^\times$, and a linear form $\ell$ such that $M . Q = Q'$, and $M . E = \alpha E'+\ell Q'$.
    Applying Lemma~\ref{prop:trick_rank4} gives $M . \tilde{E} = \tilde{M . E} = \tilde{\alpha E' + \ell Q'} = \alpha \tilde{E}'$ which proves one direction. Conversely, $(Q,\tilde{E})$ and $(Q,E)$ define the same curve, as do $(Q',\tilde{E}')$ and $(Q',E')$.
\end{proof}

As we have seen, we know how to compute linear transformations between generic cubic surfaces using covariants. Proposition~\ref{prop:isonorm} shows that we can compute the isomorphisms between the curves $C$ and $C'$ by first computing the isomorphisms between $\tilde{E}$ and $\tilde{E}'$ and then retaining the ones mapping $Q$ onto $Q'$.

This approach fails when the cubic surface defined by $\tilde{E}$ admits non-trivial automorphisms (which is always the case if $C$ does). In that case, one instead solves a system of polynomial equations to determine the coefficients of the transformation, imposing compatibility conditions on the images of $Q$, $\tilde{E}$, and the order $1$ covariants of $\tilde{E}$.
The function \texttt{IsIsomorphicGenus4()} implements this strategy, combining the covariant-based method with a Gröbner basis approach when necessary.



\medskip

For curves supported on a quadric $Q$ of rank $3$, the situation is more subtle.
We introduce two linear forms derived from $Q$ and $E$ that transform covariantly and contravariantly, respectively, under the action of $\mathrm{GL}_4(K)$. These forms are invariant under the transformation $E \mapsto E + \ell Q$.
By normalizing these quantities and applying a similar argument as in the rank $4$ case, we eliminate the residual action $E \mapsto E + \ell Q$, and apply a strategy similar to the one presented before. 

Let $\mathrm{H}(Q)$ denote the Hessian matrix of $Q$. By definition, the rank of
the quadric surface is the rank of $\mathrm{H}(Q)$.

\begin{lemma}
    Let $Q \in K[x_0,x_1,x_2,x_3]_2$ be of rank $3$, and let $v_Q$ be a nonzero element of the right kernel of $\mathrm{H}(Q)$ canonically identified to an element in $V$.
    Then for any $M\in \GL_4(K)$, $(M^{-1})^t . v_Q$ lies in the right kernel of $\mathrm{H}(M\cdot Q)$.
\end{lemma}

\begin{proof}
A direct calculation gives $\mathrm{H}(M . Q) = M^{-1}\mathrm{H}(Q)\,(M^{-1})^t$ and $(M^{-1})^t . v_Q = M^t v_Q$. The lemma follows.
\end{proof}

\begin{remark}
   The vector $v_Q$ can be interpreted as a contravariant of order $1$ of $\mathrm{Sym}^2(K^4)$, even though it is defined only up to multiplication by a scalar.
\end{remark}

Finding a transformation that sets $v_Q = v_3$ forces $Q$ to have no term in $x_3$, so that $Q\in K[x_0,x_1,x_2]_2$.
From now on, we assume $Q$ has been normalized in this way and that $v_Q = v_3$.

Now, let $E\in \mathrm{Sym}^3(V^*)$ be such that the intersection of the quadric defined by $Q$ and the cubic defined by $E$ is a smooth non-hyperelliptic curve $C$ of genus $4$.
Since $C$ is smooth, the coefficient of $x_3^3$ in $E$ is nonzero. A
direct computation also shows that the coefficient of $x_3$ in $(Q, Q, Q,
E)_2$ is nonzero. Moreover, $(Q,Q,Q,E)_2$ is a $\mathrm{GL}_4$-covariant of
$\mathrm{Sym}^2(K^4)\times \mathrm{Sym}^3(K^4)$, and it is invariant under the
action $E\mapsto E + \ell Q$. Compare this, where $(Q,Q,Q,E+\ell Q)_2 = (Q,Q,Q,E)_2$, to the rank $4$ case, where
\begin{displaymath}
(Q,Q,Q,E+\ell Q)_2 = (Q, Q, Q, E)_2 + \frac{2}{3} (Q, Q, Q, Q)_2\ell\,.
\end{displaymath}
Transforming $(Q,Q,Q,E)_2$ into $\alpha x_3$ for some $\alpha\in K^\times$ amounts to using a matrix of the form
\begin{displaymath}
  M = \begin{pmatrix}
    1 & 0 & 0 & 0\\
    0 & 1 & 0 & 0\\
    0 & 0 & 1 & 0\\
    \star & \star & \star & 1
  \end{pmatrix}\,,\quad {}^tM^{-1} = \mathrm{Com}(M) = \begin{pmatrix}
    1 & 0 & 0 & \star\\
    0 & 1 & 0 & \star\\
    0 & 0 & 1 & \star\\
    0 & 0 & 0 & 1
  \end{pmatrix}\,,
\end{displaymath}
so that ${}^tM^{-1} .  v_Q = v_3$ remains unchanged.

After dividing $E$ by the coefficient of $x_3^3$, we can write $E = x_3^3+x_3 f_2(x_0,x_1,x_2)+f_3(x_0,x_1,x_2)$.
To eliminate the action $E\mapsto E + \ell Q$, we define
\begin{displaymath}
  \tilde{f_2} = (Q,Q,Q)_2 f_2 - \frac{3}{5} (Q,Q,f_2) Q, \quad
  \tilde{f_3} = (Q,Q,Q)_2 f_3 - \frac{3}{5} (Q,Q,f_3) Q,
\end{displaymath}
where $Q$ is now seen as a polynomial in $K[x_0,x_1,x_2]$, and set $\tilde{E} = (Q,Q,Q)_2x_3^3+x_3\tilde{f_2} + \tilde{f_3}$.
Isomorphisms between two genus $4$ curves in this normal form are characterized by the following proposition.
\begin{proposition}\label{prop:isos_g4_gl3}
    Let $Q,\, Q'\in K[x_0,x_1,x_2]_2$ of rank $3$, and let $\tilde{E},\,\tilde{E'}$ be in the normal form described above, defining non-hyperelliptic curves $C$ and $C'$ of genus $4$.
Then $C$ and $C'$ are isomorphic if and only if there exist $M\in \mathrm{GL}_3(K)$ and $\alpha,\lambda,\mu\in K^\times$ such that the diagonal by block matrix $M'\in\mathrm{GL}_4(K)$ given by $\mathrm{diag}(M,\alpha)$ satisfies $M' .\tilde{E} = \lambda\tilde{E'}$ and $M' . Q = \mu Q'$.
\end{proposition}

Now that $\tilde{E}$ and $\tilde{E'}$ are normalized for the residual action $\tilde{E}\rightarrow \tilde{E}+\ell Q$, we reduce to the previous case and compute $4$ covariants of order $1$ of $\tilde{E}$.

While the covariant-based step may fail in non-generic cases, Proposition~\ref{prop:isos_g4_gl3} remains valid in full generality. In that case, one instead solves a system of polynomial equations to determine the coefficients of the transformation, imposing compatibility conditions on the images of $Q$, $\tilde{E}$, and the order $1$ covariants of $\tilde{E}$. A call to the function \texttt{IsIsomorphicGenus4()} applies the generic covariant strategy first, and resorts to a Gröbner basis when necessary.

It should in principle be possible to use the output of the previous function to compute twists via \texttt{Twist(C,H)}.

\bibliographystyle{alphaabbr}

\bibliography{synthbib}

\appendix

\section{Dixmier--Ohno invariants in positive characteristic}
\label{sec:DOpo}

    This appendix is drawn from the first author's thesis~\cite[Chap.~5]{Bou-thesis}. We refer the reader to \cite{DerKem} for  background on representation theory, in particular the theory of roots and weights. We first briefly present the theory of good filtrations, tailored to our needs.

    \subsection{Good filtrations}
        Let $n$, $d>1$, and let $p$ be a prime number. Let $k$ be an algebraically closed field of characteristic $p$.
        We are mainly interested in the  invariant algebra of $V_{n,d}(k)=\mathrm{Sym}^d(k^n)$, acted upon in the usual way by $\Gamma_k = \mathrm{SL}_n(k)$. Most, if not all, of the properties listed below extend to connected reductive groups.
        Let $T$ be the maximal torus of $\Gamma_k$ consisting of diagonal matrices of determinant $1$, and let $B$ be the Borel subgroup of $\Gamma_k$ whose elements are the upper triangular matrices with determinant $1$.
        Let $U\subseteq B$ be the unipotent subgroup of $\Gamma_k$ composed of upper triangular matrices with $1$'s on the diagonal.
        Let $\Lambda^+$ denote the set of dominant weights associated to $U$.

        \begin{definition}
            Let $\lambda\in \Lambda^+$. A \textit{dual Weyl module} $\nabla(\lambda)$ is an irreducible rational representation of $\Gamma_K$ with highest weight $\lambda$.
        \end{definition}

        \begin{remark}\label{rem:schur_irr}
        We asume that $p>d$.
        It is shown in \cite[Cor. 2.5]{DM21} (as a direct consequence of \cite[Theorems 2 and 3]{Totaro}) that for all $\lambda \vdash d$, the Schur modules $S_\lambda(k^n)$ are irreducible representations of $\mathrm{GL}_n$. Therefore, when $p>d$, Schur modules are  dual Weyl modules.
        \end{remark}
        \begin{definition}[\protect{\cite[Def. 4]{DM21}}]
            A rational representation $V$ of $\Gamma_k$ is said to have a \textit{good filtration} if there exist finite-dimensional subrepresentations $(V_i)_i$ of $V$ such that  \[0 = V_0\subseteq V_1 \subseteq V_2 \subseteq \ldots \quad\text{ and }V = \bigcup_iV_i\,,\] and each subquotient $V_i/V_{i-1}$ is a dual Weyl module.
        \end{definition}

        \begin{remark}[\protect{\cite[Rem. 5.2]{DM21}}]
            When the group $\Gamma$ can be defined over $\mathbb{Z}$, dual Weyl modules can also be defined over $\mathbb{Z}$, as $\mathbb{Z}$-modules (instead of vector spaces).
            For instance, the group $\Gamma = \mathrm{SL}_n$ can be defined over $\mathbb{Z}$, and for every algebraically closed field $k$ we have an isomorphism of algebraic groups
            \[\mathrm{SL}_n(k) \simeq \Gamma\otimes_{\mathbb{Z}}k.\]
            The dual Weyl modules can also be defined over $\mathbb{Z}$: for any $\lambda\in \Lambda^+$, there exists a free $\Gamma$-module $\nabla_{\mathbb{Z}}(\lambda)$ such that for every algebraically closed field $k$ we have an isomorphism of $\mathrm{SL}_n(k)$-representations
            \[\nabla_{\mathbb{Z}}(\lambda) \otimes_{\mathbb{Z}}k \simeq \nabla_k (\lambda).\]
            For $\Gamma = \mathrm{SL}_n$, these are the Schur modules defined over $\mathbb{Z}$.
        \end{remark}

        \begin{lemma}[\protect{\cite[Cor. 22]{DM21}}]\label{lem:good_mod}
            Let $V$ be a rational representation of a connected reductive group over $k$, such that $V$ is of dimension $n$ and admits a good filtration.
            If $p = \mathrm{char}(k)>n$, then $k[V]$ admits a good filtration.
        \end{lemma}

        \begin{corollary}\label{cor:good_filtration}
            If $p > \binom{n+d-1}{d}$, then $k[V_{n,d}(k)]$ admits a good filtration.
        \end{corollary}

        \begin{proof}
            Let $\lambda_1$ be the fundamental weight associated to the simple root \[\varepsilon_1-\varepsilon_2:\mathrm{Diag}(t_1,\ldots,t_n)\in\mathrm{SL}_n(k)\longmapsto t_1/t_2.\]
            We know that $V_{n,d}(k) = S_{(d)}(k^n)=\nabla(d\lambda_1)$, and since $p > \binom{n+d-1}{d} > d$, this Schur module is irreducible by Remark \ref{rem:schur_irr}. Thus $V_{n,d}(k)$ is a dual Weyl module.
Moreover, $p > \binom{n+d-1}{d} = \dim(V_{n,d}(k))$, therefore by Lemma~\ref{lem:good_mod}, $k[V_{n,d}(k)]$ admits a good filtration.
        \end{proof}

        Modules with good filtrations have many interesting properties, and we refer the reader to~\cite{jantzen} for further details.
        We focus here on two consequences: one concerning Cohen-Macaulayness and the other the Hilbert series.

        \begin{proposition}[\protect{\cite[Thm. 6]{hashimoto01}}]
            Let $k$ be an algebraically closed field of characteristic $p>0$, and $\Gamma$ a connected reductive group over $k$.
            Let $V$ be a finite-dimensional rational representation of $\Gamma$, and set $S = k[V]$. If $S$ has a good filtration,
            then $k[V]^\Gamma$ is Cohen-Macaulay.
        \end{proposition}

        In particular, in characteristic greater than the dimension of the rational representation considered, the algebra of invariants is Cohen–Macaulay, just as in characteristic $0$. We obtain a similar result for the Hilbert series.

        \begin{proposition}[\protect{\cite[Prop. 8.7]{DM21}}]\label{prop:hilb_ser_pos_char}
            Let $V$ be a rational representation of dimension $n$ of a connected reductive group $\Gamma$ defined over $\mathbb{Z}$ which admits a good filtration.
            Let $k$ be an algebraically closed field such that $k[V]$ admits a good filtration.
            Then the Hilbert series of $k[V]^\Gamma$ is the same as the Hilbert series of $\mathbb{C}[V]^\Gamma$.
        \end{proposition}

        \begin{corollary}\label{cor:general}
            If $p > \binom{n+d-1}{d}$, then \[\mathrm{Hilb}\left(k[V_{n,d}(k)]^{\mathrm{SL}_n(k)},t\right) = \mathrm{Hilb}\left(\mathbb{C}[V_{n,d}(\mathbb{C})]^{\mathrm{SL}_n(\mathbb{C})},t\right)\,.\]
        \end{corollary}

        \begin{proof}
            We follow~\cite[Sec. 5]{DM21}. The algebraic group $\mathrm{SL}_n$ is connected, reductive, and defined over $\mathbb{Z}$.
            The representation $V_{n,d}$ is a dual Weyl module, hence admits a good filtration.
            If $p > \binom{n+d-1}{d}$, then $k[V_{n,d}(k)]$ admits a good filtration.
            Proposition~\ref{prop:hilb_ser_pos_char} then yields the result.
        \end{proof}

        \subsection{Application to the invariant ring of ternary quartics}

        \begin{proposition}[\cite{shioda67}]\label{prop:hilbert-ser}
            The Hilbert series of the algebra of invariants of ternary quartics over an algebraically closed field of characteristic $0$ is \[\mathrm{Hilb}\left(\mathbb{C}[V_{3,4}(\mathbb{C})]^{\mathrm{SL}_3(\mathbb{C})}, t\right) = \frac{P(t)}{(1-t^3)(1-t^6)(1-t^9)(1-t^{12})(1-t^{15})(1-t^{18})(1-t^{27})}\,,\]
            where
            \begin{align*}
            P(t)&=1+t^9+t^{12}+t^{15}+2t^{18}+3t^{21}+2t^{24}+3t^{27}+4t^{30}+
            3t^{33}+4t^{36}+4t^{39}+3t^{42}+4t^{45}+3t^{48}\\
                &\quad +2t^{51}+3t^{54}+2t^{57}+t^{60}+t^{63}+t^{66}+t^{75}.
            \end{align*}
        \end{proposition}

        \begin{theorem}\label{thm:DO_car_pos}
            Let $k$ be an algebraically closed field of characteristic $p$.
            When $p > 13$, the algebra of invariants $k[V_{3,4}(k)]^{\mathrm{SL}_3(k)}$ is generated by the Dixmier--Ohno invariants.
        \end{theorem}

        \begin{proof}
            By Corollary~\ref{cor:general}, we have \[\mathrm{Hilb}\left(k[V_{3,4}(k)]^{\mathrm{SL}_3(k)}, t\right) = \mathrm{Hilb}\left(\mathbb{C}[V_{3,4}(\mathbb{C})]^{\mathrm{SL}_3(\mathbb{C})}, t\right).\]

            Furthermore, a primary system of invariants of $k[V_{3,4}]^{\mathrm{SL}_3(k)}$ was established for every $p$ greater than $7$ in~\cite{LLLR21}.
            For $p\notin \{19, 47, 277, 523\}$, the Dixmier invariants still form a primary system of invariants.
            For the exceptional primes $p\in \{19, 47, 277, 523\}$, a different homogeneous system of parameters (of the same degrees) is obtained by replacing the degree $9$ invariant $I_9$ with $I_9+J_9$, where $J_9$ is the degree $9$ Ohno invariant.

            Since for each characteristic $p>7$ we know a primary system of invariants of degrees matching those of the denominator of the Hilbert series, checking that the Dixmier--Ohno invariants $\mathrm{DO}$ generate the algebra of invariants reduces to
            verifying that for every degree $d \leq 75$ (the degree of the numerator of the Hilbert series), we have
            \[\mathrm{dim}(k[\mathrm{DO}]_{d}) = \mathrm{dim}\left(k[V_{3,4}(k)]^{\mathrm{SL}_3(k)}_{d}\right)\,.\]
            Since it is impossible to check this property directly for every $p > 13$, we perform the verification over $\mathbb{Z}\left[\frac{1}{13!}\right]$.
            Specifically, we evaluate a conjectured basis of invariants for each degree $d\leq 75$ on a set of ternary quartics with small integer coefficients, in a number exceeding the expected dimension.
            Computing the Smith normal form of the resulting matrix then
            identifies, at least, the primes for which the rank is not what is expected.

            For larger degrees, computing the full Smith normal form is computationally intensive, so we compute the greatest common divisor of several maximal minors of the matrix instead. This gcd is a multiple of the last diagonal element of the Smith normal form, and hence also  of every diagonal element. Computing the gcd of three maximal minors reveals that the divisors  are contained in  $\{2,3,5,7\}$ for all degree up to $75$.

            For $p\in \{19, 47, 277, 523\}$, one checks separately that the algebra generated by the Dixmier--Ohno invariants has the correct dimension for all degrees between $0$ and $75$, using the Hilbert series and a similar  evaluation/rank strategy, but this time over a finite field (which is much faster).

            We  conclude that in characteristic  $p$ greater than $13$, the dimension of the homogeneous components of $k[V_{3,4}(k)]^{\mathrm{SL}_3(k)}$ coincides with that of the algebra generated by the Dixmier--Ohno invariants.
    By inclusion , this implies equality.
        \end{proof}

        \begin{remark}
            We proved that for every algebraically closed field $k$ of characteristic $p> 13$,  \[k[V_{3,4}(k)]^{\mathrm{SL}_3(k)} = \mathbb{Z}[\mathrm{DO}]\otimes_{\mathbb{Z}} k.\]
            Since we only observed possible dimension drops for $p\leq 7$, one may hope that the Hilbert series remains unchanged for $p = 11$ or $13$, in which case  Theorem~\ref{thm:DO_car_pos} would extend immediately to $p > 7$.
        \end{remark}

        This observation naturally leads to the question of whether the bound $p > \binom{n+d-1}{d}$ is sharp, a question we hope will receive a negative answer in the future.
        \medskip

        \begin{question}
            Is the bound $p > \binom{n+d-1}{d}$ in Corollary~\ref{cor:good_filtration} sharp for ensuring that $k[V_{n,d}(k)]$ admits a good filtration?
        \end{question}
\end{document}